\documentclass[letterpaper, 10 pt, conference]{ieeeconf}  

\IEEEoverridecommandlockouts                              
\usepackage[dvipdfm]{graphicx} 
\usepackage{subfigure}
\usepackage{amsmath}
\usepackage{amsfonts}
\usepackage[psamsfonts]{amssymb}
\usepackage{color}
\usepackage{comment}
\usepackage{times} 

\newtheorem{proposition}{Proposition}

\newtheorem{assumption}{Assumption}

\usepackage{epsfig} 

\usepackage{amsmath,amssymb,graphics,epsfig,subfigure}

\newcommand{\be}{\begin{equation}}
\newcommand{\ee}{\end{equation}}
\newcommand{\bea}{\begin{eqnarray}}
\newcommand{\eea}{\end{eqnarray}}
\newcommand{\ba}{\begin{array}}
\newcommand{\ea}{\end{array}}
\newcommand{\beas}{\begin{eqnarray*}}
\newcommand{\eeas}{\end{eqnarray*}}
\newcommand{\leftm}{\left[\begin{array}}
\newcommand{\rightm}{\end{array}\right]}

\title{\LARGE \bf
Infinitesimal Perturbation Analysis of Stochastic Hybrid Systems: Application to  Congestion Management in  Traffic-Light Intersections
}

\author{Y. Wardi$^*$\thanks{$^*$School of Electrical and Computer Engineering, Georgia Institute of
Technology, Atlanta, GA 30332, USA. Email: ywardi@ece.gatech.edu.}$^{\dag}$ and C. Seatzu$^{\#}$\thanks{$^{\#}$Department of Electrical and Electronic
Engineering, University of Cagliari, Italy.  Email:
seatzu@diee.unica.it.}\thanks{$^{\dag}$ Research supported in
part by  NSF under Grant CNS-1239225.}
}

\begin{document}

\maketitle
\thispagestyle{empty}
\pagestyle{empty}


\begin{abstract}
This paper presents a new approach to congestion management at traffic-light intersections. The approach is based on
controlling the relative lengths of red/green
cycles in order to have the congestion level track a given reference. It uses an integral
control with adaptive gains, designed to provide fast tracking and wide stability margins. The gains are
inverse-proportional to the derivative of the plant-function with respect to the control parameter, and
are computed by
infinitesimal perturbation analysis. Convergence of this technique is shown to be robust with respect  to modeling uncertainties, computing errors, and other random effects. The framework is presented in the setting of stochastic hybrid
systems, and applied to a particular traffic-light model. This is but an initial study and hence the latter model is simple, but
it captures some of the salient features of traffic-light processes. The paper concludes with comments on possible extensions of the proposed approach
to traffic-light grids with realistic flow models.
\end{abstract}

\section{Introduction}
Infinitesimal Perturbation Analysis (IPA) is an established sample-path technique for sensitivity estimation of performance functions
defined on discrete event dynamic systems and stochastic hybrid systems. In a typical scenario $L(\theta)$ is a real-valued, random function of a parameter $\theta\in R^n$, defined over a common probability space
$(\Omega,{\cal F},P)$, and for a particular realization corresponding to $\omega\in\Omega$,  IPA computes its sample
($\omega$-dependent) derivative (gradient) $\nabla L(\theta)$. This sample gradient can act, under certain
circumstances,  as an estimator of the gradient of the expected-value function $J(\theta):=E\big[L(\theta)\big]$, with
$E[\cdot]$ denoting expectation in $(\Omega,{\cal F},P)$. This can be used, in conjunction with gradient-descent
algorithms such as stochastic approximation, to minimize the function
$J(\theta)$ to the extent of computing a local minimum. For extensive presentation of IPA and its scope
in optimization, , please see
\cite{Ho91,Glasserman91,Cassandras99}.

In recent years there has been a mounting interest in the application of IPA to stochastic hybrid systems, and
especially to stochastic flow networks, comprising generalizations of fluid queues \cite{Cassandras02,Cassandras06}. The reason is that, for an extensive
class of performance functions in this setting, IPA was shown to be computable via simple algorithms directly from quantities that
are observable from realizations of the state of the system. Moreover, the resulting gradient estimators are quite
robust to modeling variations. For these reasons, for a variety of
applications, the IPA gradients arguably can be computed in real time without having concerns about the accuracy of the
underlying models. References \cite{Cassandras10,Wardi10} provide detailed
discussions of these points as well as unified frameworks
for IPA in the setting of stochastic hybrid systems.
Until recently the main areas mentioned for potential applications of IPA were in manufacturing and telecommunications,
but  lately
 there has been a growing interest in transportation networks as well, and especially in traffic-light control (see \cite{Fleck14} for
 a survey).

 The main objective
 of traffic-light control is to reduce, or minimize congestion at
 traffic intersections. Early techniques developed for these control and optimization problems include
 dynamic programming \cite{Porche96} and   linear complementary algorithms \cite{DeSchutter99},
 while
 more recent approaches are based on  Markov-decision processes \cite{Yu06}, game theory \cite{Alvarez10},
 and
 mixed-integer programming \cite{Dujardin11};
 \cite{Papageorgiou03} contains an early survey, and a recent one can be found in  \cite{Fleck14}.   Regarding applications of IPA,
 early results were presented in
\cite{Fu03,Panayiotou05}, and a recent systematic approach has been
developed in
 \cite{Geng12a,Geng12b,Geng13,Fleck14}. This approach defines the traffic-light control problem in the aforementioned
 setting of stochastic hybrid systems \cite{Cassandras10}, and develops for it effective  IPA-based algorithms.

The development of IPA since its inception has been motivated primarily by applications to
performance optimization in discrete event and hybrid systems. This paper follows a different track in pursuing
an application to  performance regulation. The term ``regulation'' here means real-time tracking of
a set (reference) performance index by tuning a control parameter, and it is a common engineering practice. In
particular, following a system's  design or optimization  with an imprecise system model, regulation
can be used to ensure that performance meets  specifications under changing system characteristics and operating
environments.

 In devising our regulation technique we aim at effective and efficient real-time implementation. Effectiveness means that the
set-point tracking algorithm is to have fast convergence under a wide set of system parameters, efficiency means simple implementation
 requiring low computing efforts, and the real-time requirement means that all input parameters to the controller be measurable
 by observing the system's state. High degree of efficiency means that we may have to tilt the balance between  speed and
 precision
 of computation in favor of the former requirement. Consequently  we design the controller for maximum speed
 and simplicity of computations, possibly by using imprecise models,
 while guaranteeing its robustness under large variations in the system's parameters.  {\it The main
 contribution of this paper is in a control system with all of these properties. It is based on an integral control
 with a variable gain, computed by using the IPA derivative of the plant function with respect to the control variable.}

 The integral controller is explained in Section II in an abstract setting.  Its application to a relevant example
 of traffic-light control, including the derivation of the IPA derivative, is presented in Section III. Section IV
 contains simulation results, and Section V concludes the paper.

 \section{Regulation Algorithm: Integral Control with Adaptive Gain}
 Consider the single-input-single-output discrete-time
 control system shown in Figure~1, where $r$ is the reference (set point) input, $y_{n}$ is the output, $e_{n}$ is the error signal, and $u_{n}$ is the input to the plant. Suppose that the plant is a time-varying, memoryless nonlinearity
 of the form
 \begin{equation}
 y_{n}=g_{n}(u_{n}),
 \end{equation}
 where $g_{n}:R\rightarrow R$, $n=1,2,\ldots$, is called the {\it plant function}.
 Given a reference input $r$, the purpose of the control system is to ensure that
 $\lim_{n\rightarrow\infty}y_{n}=r$. To this end it is natural to choose the controller to be an integrator having,
 for example, the transfer function $G_{c}(z)=Az^{-1}/(1-z^{-1})$ for a given gain $A>0$.
 However, integral controllers may display oscillatory behavior and have narrow stability margins. Furthermore,
 due to the time-varying nature of the system it may not be easy (or possible) to choose a gain that fits all
 possible scenarios.
 For this reason we use an integral control with adaptive gain, $A_{n}$, having the time-domain representation
 \begin{equation}
 u_{n}\ =\ u_{n-1}+A_{n}e_{n-1},
 \end{equation}
where the error signal is defined as
 \begin{equation}
 e_{n}=r-y_{n}.
 \end{equation}
 We choose the gain $A_{n}$ to be defined via the equation
 \begin{equation}
 A_{n}\ =\ \frac{1}{g_{n}^{\prime}(u_{n-1})},
 \end{equation}
 with ``prime'' denoting derivative. Equations (1) -- (4), computed cyclically in the order
 $(4)\rightarrow(2)\rightarrow(1)\rightarrow(3)$, define the dynamics of the closed-loop system.
 In fact, we implicitly consider the following scenario in a real-time setting: the quantities $u_{n-1},\ y_{n-1}$, and
 $e_{n-1}$ have been computed or derived by the start of the $n$th computing cycle, and at that time the following sequence of
 operations takes place and are completed before the end of the cycle: (i) $A_{n}$ is computed by Equation (4);
 (ii) $u_{n}$ is computed by the controller via Equation (2); and (iii) the system yields  $y_{n}$ and $e_{n}$ (see (1), (3)).
 The main computation is that of $A_{n}$ in Equation (4), since an exact evaluation of
 $g_{n}^{\prime}(u_{n-1})$ may not be possible in real time. In this case it may be necessary to trade off precision with
 computational expediency, with the result that Equation (4) be replaced by
 \begin{equation}
 A_{n}\ =\ \frac{1}{g_{n}^{\prime}(u_{n-1})+\xi_{n}},
 \end{equation}
 where $\xi_{n}$ is a  quantity representing the computing error. We will say more about this point in the sequel.

 \begin{figure}[!t]
\centering
\includegraphics[width=3in]{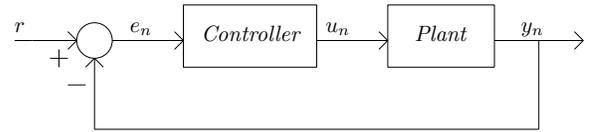}
\caption{Control System}
\label{Figure 1}
\end{figure}

 The rationale behind Equation (4) can be seen by considering the time invariant plant, where $g$ is not a function
 of $n$ and hence Equation (1) has the form $y_{n}=g(u_{n})$. In this case, it can be seen that
 the control loop, defined via Equations (1) --(4), implements Newton's method for solving the equation
 $g(u)=r$, known to converge, in the sense that
 \begin{equation}
 \lim_{n\rightarrow\infty}e_{n}\ =\ 0,
 \end{equation}
 under broad assumptions. In the event that Equation (5) is used instead of (4),
 namely  (still assuming a time-invariant plant) $A_{n}=1/(g^{\prime}(u_{n-1})+\xi_{n})$,
 the regulation scheme still converges in the sense of Equation (6),  under broad assumptions,
 as long as the relative error is under 100\%. For example, if the function $g$ is concave or convex,
 and if, for some $\alpha\in(0,1)$, the relative error of the plant derivative, defined as $\varepsilon_{n}:=|\xi_{n}|/|g_{n}^{\prime}(u_{n+1})|$,
 satisfies the inequality
 $\varepsilon_{n}<\alpha$ for all $n=1,2,\ldots$, then Equation (6) is satisfied (see \cite{Almoosa12}).

 In the case of a time-varying plant defined by
 Equation (1),
 Equation (6)  no longer can be expected. However, the error sequence $\{e_{n}\}$
 asymptotically gets close to 0 by an amount that depends, monotonically, on a measure of the
 system's variability. For example, \cite{Almoosa12} proved, under the convexity assumption, that for every
 $\epsilon>0$ there exists $\delta>0$ such that, if $|g_{n}(u_{n-1})-g_{n-1}(u_{n-1})|<\delta$ for every $n=1,2,\ldots$, then
 \begin{equation}
 \limsup_{n\rightarrow\infty}|e_{n}|\ <\ \epsilon.
 \end{equation}
 We point out that the convexity assumption can be relaxed to local convexity or concavity as long as guards are put in
 place to ensure an appropriate upper bound on the terms $|u_{n}-u_{n-1}|$. Moreover, bounds on the relative errors
 $\varepsilon_{n}:=|\xi_{n}|/|g_{n}^{\prime}(u_{n-1})|$ practically need not be computed a priori but can be verified  from system
 simulation as  done in \cite{Almoosa12a}.

 This control law was applied to regulate the dynamic core-power in computer processors  by the applied frequency.
 The plant, comprising the  frequency-to-power relationship, has an established model based on
 physical principles \cite{Floyd07}, represented by
  a  convex, time-varying, memoryless
 nonlinearity. Its time variability is due to the activity factor of the program load, a quantity representing the amount
 of  switching activity of the logic gates at the core. This quantity generally is unpredictable and cannot be measured
 in real time, and hence the plant functions cannot be computed. Although the output $y_{n}$ can be measured, the derivative $g_{n}^{\prime}(u_{n-1})$ requires a formula.
 It turns out that this derivative is computable from quantities that can be measured in real time and hence
 the regulation scheme could be applied; for simulation results with an industry-grade simulator,
 please see \cite{Almoosa12}.

 Reference \cite{Almoosa12a} considers regulating the instruction throughput in a similar core, also as a function
 of frequency. Each control cycle consists of a fixed number of clock cycles and takes about 10 miliseconds, during which
 the applied frequency is fixed. The output $y_{n}$ is defined as the average throughput over a given cycle,
 and the model for the plant is a queueing network representing instruction-processing at the core.
 The queueing model has multiple precedence constraints and is complicated in various other ways,
  and hence defies analysis for deriving closed-form formulas for the plant function. However, we estimated the derivative term
  $g_{n}^{\prime}(u_{n-1})$ by using IPA, which yielded simple formulas that could be computed in real time. We point
  out that these IPA derivatives were statistically biased and hence ``wrong'', but extensive simulations
  on various benchmark systems yielded maximum relative error of 30\%. Encouraged by these findings we applied the regulation scheme despite the bias of IPA, and the results can be seen in \cite{Almoosa12a}.

  This example highlights  two  salient points of our proposed regulation framework. First, the plant is modeled as a queueing
  network which is a highly dynamic system, but the cost function,   consisting of the average throughput over
  a certain amount of time, allows us to consider it as a {\it memoryless, time-varying} nonlinearity. The time variability
  is due to the uncertain, random element as well as to variations  that are inherent in the system's dynamics.
  This view serves us better than the dynamic view by dint of Equation (7), and suggests that, under conditions of stochastic stability,  longer control cycles
should result in smaller values of $\lim_{n\rightarrow\infty}|e_{n}|$.
Second, the important feature of IPA and its use in Equation (4) is its simplicity and on-line computability, and therefore
the main objective of the analysis below is to derive simple terms for the sample derivatives
$g_{n}^{\prime}(u_{n-1})$.

\section{Traffic-Light Regulation Problem}
This section concerns  an application of the  aforementioned control technique to a traffic-light intersection model.
In order to highlight the salient features of the regulation scheme, we consider only a simple model and defer
 a discussion  of more detailed models to a future study. Thus, consider a traffic-light  intersection of two unidirectional (one-way) roads. Suppose
that the light in each direction alternates between red and green signals,
and  for the sake of simplicity we assume that there is no orange light,
and hence,  epochs of red signal in one direction (road) correspond to green signal in the other direction.
Let us focus attention on one of the roads  and define a {\it light cycle} as a red period
followed by a green period (for the other road, the same cycle is comprised of green followed by red).
Suppose that the light cycle time is fixed, and denoted by $C$, and let the length of the red period comprise
the control variable, $\theta\in[0,C]$.

One of the ways to characterize congestion is by the traffic buildup in front of a traffic  light,
and as in \cite{Geng12a,Geng13,Fleck14}, we model its dynamics  by a fluid queue. Such a queue is driven by two random processes: the arrival rate
and the service rate. The arrival-rate process does not depend on $\theta$ and hence it  is denoted by
$\{\alpha(t)\}$, while the service-rate process depends on $\theta$ in a manner described below, and hence it
is denoted by $\{\beta(\theta;t)\}$. Given  an integer $N>0$ we  define the performance function that we seek to regulate
as  the time-average of the buffer contents (amount of fluid at the buffer) during $N$ light cycles. Defining $T:=NC$, we  call the period $[0,T]$ a {\it control cycle}, and we note that a control cycle consists of $N$ light cycles.
We assume that $\theta$ remains fixed during each
control cycle and it is changed, by the regulation process,
only at the boundary points between  consecutive control cycles.

Consider  a control cycle at a given $\theta\in[0,C]$. We define the  buffer contents during the cycle,
denoted by $x(\theta,t)$, by the one-sided differential equation
\begin{equation}
\frac{dx}{dt^{+}}(\theta;t)=\left\{
\begin{array}{ll}
\alpha(t)-\beta(\theta;t), & {\rm if}\ x(\theta,t)>0\\
0, & {\rm if}\ x(\theta,t)=0,
\end{array}
\right.
\end{equation}
where, for the sake of simplicity, we assume the initial condition $x(\theta,t)=0$
(our simulations use a different initial condition as will be explained in the sequel).
We assume that, $\forall\ \theta\in[0,C]$, w.p.1, $\alpha(t)$ and $\beta(\theta,t)$ are piecewise continuously
differentiable (but not necessarily continuous) in the interval
$t\in[0,T]$, and this ensures that
 $x(\theta,t)$ is well defined by Equation (8).
 The performance function $L(\theta)$ is defined as
 \begin{equation}
 L(\theta)=\frac{1}{T}\int_{0}^{T}x(\theta;t)dt.
 \end{equation}
 Given a reference set-point $r>0$, the objective of the regulation scheme is to compute a sequence
 of control variables $\theta_{1},\theta_{2},\ldots$, such that the sequence $L(\theta_{1}),L(\theta_{2}),\ldots$ tracks
 $r$ as best as possible.

Regarding the service rate process, we define a model that includes the cases where, upon the light switching from
red to green, the service rate either jumps to, or ramps up towards a given maximum value.
In the latter case  the ramp-up period  depends  on the
queue length at the time the green epoch starts,  and lasts until either one of the following two
events occur: (i) The buffer becomes empty, or (ii) the light switches back to red.
In the first event the service rate jumps to the maximum value, and in the second event, it jumps down to 0.
We model the ramp-up rate process as a random function to account for fluctuations that are hard
to model, and assume mutually independent realizations of it in
successive cycles.

Formally, let $\beta_{m}>0$ be a given constant, and let  $b(t):[0,C]\rightarrow[0,\beta_{m}]$
be a monotone-nondecreasing random function that does not depend on $\theta$.
 On a given light cycle
$[kC,(k+1)C)$, we define the service rate of the queue as follows,
\begin{equation}
\beta(\theta,t)=\left\{
\begin{array}{ll}
0, & {\rm if}\ t\in[kC,kC+\theta)\\
b(t-(kC+\theta)), & {\rm if}\ t\in[kC+\theta,kC+C),\\
& {\rm and}\ x(\theta,t)>0\\
\beta_{m}, & {\rm if}\ t\in[kC+\theta,kC+C),\\
 &  {\rm and}\ x(\theta,t)=0.
\end{array}
\right.
\end{equation}
Note that this includes the case where $\beta(\theta,t)$ jumps directly from $0$ to $\beta_{m}$ when the light turns green.

In the rest of this section we derive the IPA formula for the derivative term $L^{\prime}(\theta,t)$. This will be
done under the following assumption.
\begin{assumption}
(i). The random functions $\alpha(t)$ and $b(t)$ are independent of each other.
(ii). W.p.1,  $\alpha(t)$ is  piecewise monotone (nondecreasing/non-increasing) and piecewise continuously differentiable in
 $t$.
 (iii). W.p.1,  the function $b(t)$ is monotone nondecreasing  and piecewise continuously differentiable in $t$.
 (iv). For every $\theta\in[0,T]$, w.p.1 $\alpha(t)$ is continuous at the points $kC$ and $kC+\theta$, $k=0,1,\ldots$.
 (v). For every $\theta\in[0,T]$, w.p.1 the function $\alpha(t)$ is continuous at any point where $\beta(\theta;t)$ is discontinuous.
 (vi). For every $\theta\in[0,C)$, w.p.1,  for every open interval $I\subset[0,T]$, it is impossible to have the relation
 $\alpha(t)=\beta(\theta,t)$ $\forall\ t\in I$ except for the case where $\alpha(t)=\beta(\theta;t)=0$.
 \end{assumption}

 Similar assumptions are routinely made in the literature on IPA in the setting of stochastic hybrid systems;
 e.g., \cite{Cassandras02,Wardi10}. We remark that the attribute ``piecewise continuously differentiable'' means that it is continuously differentiable at all but
 a finite set of time-points $t$. At those points it may be discontinuous.
 This set of points  may depend on the sample $\omega$, and its cardinality   need not be upper-bounded over
 $\omega\in\Omega$.

In the forthcoming discussion we use the `prime' notation for derivatives with respect to
$\theta$, and the `dot' notation for derivatives with respect to $t$. Thus,
$x^{\prime}(\theta;t):=\frac{\partial x}{\partial\theta}(\theta;t)$, while $\dot{x}(\theta;t)=\frac{\partial x}{\partial t}(\theta;t)$.

Fix $\theta\in[0,T]$. $x(\theta,t)$ is continuous in $t$ by Equation (8), and by (9),
\begin{equation}
L^{\prime}(\theta)=\frac{1}{T}\int_{0}^{T}x^{\prime}(\theta;t)dt.
\end{equation}
We next derive formulas for the derivative term $x^{\prime}(\theta;t)$.

First, suppose
 that $t$ lies in the interior of an empty period
(namely, the continuous-queue analogue of an idle period; a maximal period when
$x(\theta;\cdot)=0$; see \cite{Cassandras02}). Then,
obviously, $x^{\prime}(\theta,t)=0$.

 Next, consider $t$ lying in a non-empty period
(the complementary of empty periods; a supremal interval where $x(\theta;\cdot)>0$). Let $u_{t}(\theta)$ be the
starting time of the nonempty period, namely, $u_{t}(\theta):=\max\{u\leq t:x(\theta;u)=0\}$;
if no such $u$ exists, $u_{t}(\theta):=0$. Let $\ell C,(\ell+1)C,\ldots,mC$ denote the time-points in the interval
$(u_{t}(\theta),t)$ when the light switches from green to red, for some integers $\ell\geq 1$ and $m\geq \ell$.
We have the following result.

\begin{proposition}
Fix $\theta\in(0,C)$, and consider $t\in[0,T]$ such that $\beta(\theta;\cdot)$
is continuous at $t$,
and  $x(\theta;t)>0$. Then
the term $x^{\prime}(\theta;t)$ has the following form,
\begin{equation}
x^{\prime}(\theta;t)=\sum_{k=\ell}^{m}\beta(\theta;(kC)^-)+\beta(\theta;t)-\beta(\theta;u_{t}(\theta)^+).
\end{equation}

For the proof we provide the diagram in Figure~2 as a visual aid.

\begin{figure}
\centering
\includegraphics[width=3.3in]{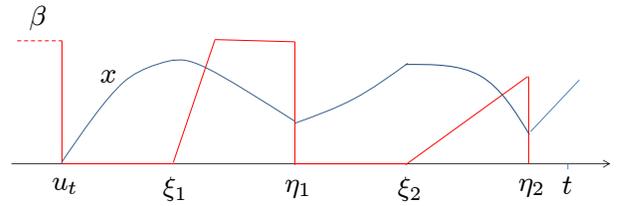}
\caption{Evolution of $\beta(\theta;\cdot)$ (red curve) and $x(\theta;\cdot)$ (blue curve):
$u_{t}=(\ell-1)C$;
$\xi_{1}=(\ell-1)C+\theta$;
 $\eta_{1}=\ell C$;
  $\xi_{2}=\ell C+\theta$;
  $\eta_{2}=(\ell+1)C$; $p=2$.
  }
\label{Figure 2}
\end{figure}

\begin{proof}
By definition of $u_{t}(\theta)$ and Equation (8), we have that
\begin{equation}
x(\theta;t)=\int_{u_{t}(\theta)}^{t}\big(\alpha(\tau)-\beta(\theta;\tau)\big)d\tau.
\end{equation}
If $\tau$ lies in the interior of a red-signal period then
$\beta(\theta;\tau)=0$ and hence $\frac{\partial}{\partial\theta}\big(\alpha(\tau)-\beta(\theta;\tau)\big)=0$;
therefore only green-signal periods need be considered in the computation of $x^{\prime}(\theta;t)$ in the following way.
Let us denote the green-signal periods in the interval $[u_{t}(\theta),t]$ by $\Gamma_{j}$, $j=1,\ldots,p$, in increasing
order, and let $\xi_{j}(\theta)$ and $\eta_{j}(\theta)$ be the boundary points of $\Gamma_{j}$, so that
$\Gamma_{j}=[\xi_{j}(\theta),\eta_{j}(\theta)]$.
Define the functions $f_{j}(\theta)$, $j=1,\ldots,p$, as follows:
\begin{equation}
f_{j}(\theta):=\int_{\xi_{j}(\theta)^{-}}^{\eta_{j}(\theta)^{+}}\big(\alpha(\tau)-\beta(\theta;\tau)\big)d\tau,
\end{equation}
with only the following  two possible exceptions:  (i) For $j=1$, if $\xi_{1}(\theta)=u_{t}(\theta)$, then
the left - limit point of the integral in (14) is $u_{t}(\theta)$ and not $u_{t}(\theta)^-$;
and (ii) for $j=p$, if $\eta_{p}(\theta)=t$ then the right-limit point of the integral is
$t$ and not $t^+$. The reason for these exceptions is that the integral in (13) is taken over $[u_{t}(\theta),t]$ and its integrant
may not be valid outside this interval.  Furthermore, it is obvious that
$x(\theta;t)=\sum_{j=1}^{p}f_{j}(\theta)$, and hence,
\begin{equation}
x^{\prime}(\theta;t)=\sum_{j=1}^{p}f_{j}^{\prime}(\theta).
\end{equation}
We next derive formulas for $f_{j}^{\prime}(\theta)$.

Consider a typical green-signal period $\Gamma:=[\xi(\theta),\eta(\theta)]\subset[u_{t}(\theta),t]$, and define
\begin{equation}
f(\theta)=\int_{\max\{\xi(\theta)^-,u_{t}(\theta)\}}^{\min\{\eta(\theta)^+,t\}}\big(\alpha(\tau)-\beta(\theta;\tau)\big)d\tau;
\end{equation}
note that this is like one of the functions $f_{j}(\theta)$ defined in (14) with the noted exceptions.
Let $\tau_{i}(\theta)$, $i=1,\ldots,q$, be the jump-times of $\beta(\theta;\cdot)$ in increasing order, in
the interval $(\xi(\theta),\eta(\theta))$. Suppose first that $\xi(\theta)>u_{t}(\theta)$ and $\eta(\theta)<t$, the cases where
$\xi(\theta)=u_{t}(\theta)$ or $\eta(\theta)=t$ will be considered later.
Taking derivative with respect to $\theta$ in (16), we obtain,
\begin{eqnarray}
f^{\prime}(\theta)\ =\
-\int_{\xi(\theta)}^{\tau_{1}(\theta)}\beta^{\prime}(\theta;\tau)d\tau\nonumber\\
-\sum_{i=1}^{q-1}\int_{\tau_{i}(\theta)}^{\tau_{i+1}(\theta)}\beta^{\prime}(\theta;\tau)d\tau
-\int_{\tau_{q}(\theta)}^{\eta(\theta)}\beta^{\prime}(\theta;\tau)d\tau\nonumber\\
-\sum_{i=1}^{q}\big(\beta(\theta;\tau_{i}(\theta)^{-})-\beta(\theta;\tau_{i}(\theta)^{+})\big)\tau_{i}^{\prime}(\theta)\nonumber\\
+\Big(\big(\alpha(\xi(\theta)^-)-\beta(\theta;\xi(\theta)^-)\big)\nonumber\\
-\big(\alpha(\xi(\theta)^+)
-\beta(\theta;\xi(\theta)^+)\big)\Big)\xi^{\prime}(\theta)\nonumber\\
+\Big(\big(\alpha(\eta(\theta)^-)-\beta(\theta;\eta(\theta)^-)\big)\nonumber\\
-\big(\alpha(\eta(\theta)^+)
-\beta(\theta;\eta(\theta)^+)\big)\Big)\eta^{\prime}(\theta);
\end{eqnarray}
we note that discontinuities in $\alpha(\cdot)$ do not matter since the process $\{\alpha(\cdot)\}$ is independent of
$\theta$, and further notice that $\alpha(\cdot)$ is continuous at jump-points of $\beta(\theta;\cdot)$ by Assumption
1.v.

Consider the integral terms in the RHS of (17). Since $\Gamma$ is a green-signal period contained in a nonempty period of
the queue, Equation (10) implies that, for every $\tau\in\Gamma$, $\beta(\theta;\tau)=b(\tau-kC-\theta)$ for some $k=1,2,\ldots$.
This implies that $\beta^{\prime}(\theta;\tau)=-\dot{\beta}(\theta;\tau)$, which allows
us to compute the integrals in the following way:
\begin{eqnarray}
-\int_{\xi(\theta)}^{\tau_{1}(\theta)}\beta^{\prime}(\theta;\tau)d\tau
=\int_{\xi(\theta)}^{\tau_{1}
(\theta)}
\dot{\beta}(\theta;\tau)\nonumber\\
=\beta(\theta;\tau_{1}(\theta)^-)-\beta(\theta;\xi(\theta)^+),
\end{eqnarray}
and similarly, for the rest of the integrals,
\begin{equation}
-\int_{\tau_{i}(\theta)}^{\tau_{i+1}(\theta)}\beta^{\prime}(\theta;\tau)=\beta(\theta;\tau_{i+1}(\theta)^-)-\beta(\theta;\tau_{i}(\theta)^+),
\end{equation}
and
\begin{equation}
-\int_{\tau_{q}(\theta)}^{\eta(\theta)}\beta^{\prime}(\theta;\tau)=\beta(\theta;\eta(\theta)^-)-\beta(\theta;\tau_{q}(\theta)^+).
\end{equation}
Substituting from Equations (18)-(20) in (17) we obtain,
\begin{eqnarray}
f^{\prime}(\theta)=\sum_{i=1}^{q}\big(1-\tau_{i}^{\prime}(\theta)\big)\big(\beta(\theta;\tau_{i}(\theta)^{-})
-\beta(\theta;\tau_{i}(\theta)^{+})\big)\nonumber\\
-\beta(\theta;\xi(\theta)^+)+\beta(\theta;\eta(\theta)^-)\nonumber\\
+\Big(\big(\alpha(\xi(\theta)^-)-\beta(\theta;\xi(\theta)^-)\big)\nonumber\\
-\big(\alpha(\xi(\theta)^+)
-\beta(\theta;\xi(\theta)^+)\big)\Big)\xi^{\prime}(\theta)\nonumber\\
+\Big(\big(\alpha(\eta(\theta)^-)-\beta(\theta;\eta(\theta)^-)\big)\nonumber\\
-\big(\alpha(\eta(\theta)^+)
-\beta(\theta;\eta(\theta)^+)\big)\Big)\eta^{\prime}(\theta).
\end{eqnarray}
Next, we observe that a green-signal period contained in the interval $(u_{t}(\theta),t)$ ends  at a time-point
$kC$, $k=1,\ldots$, meaning that $\eta(\theta)=kC$ which is independent of $\theta$, and hence
$\eta^{\prime}(\theta)=0$. This implies that the last additive term in Equation (21) is zero.
Furthermore, each time $\tau_{i}(\theta)$ lies in the interior of the green-signal period $\Gamma$, and hence there exists an open interval
containing it where, by Equation (10), $\beta(\theta,\tau)=b(\tau-kC-\theta)$; this implies that $\tau_{i}^{\prime}(\theta)=1$, which annuls
the first additive term in the RHS of (21). All of this reduces (21) to the following equation,
\begin{eqnarray}
f^{\prime}(\theta)=
-\beta(\theta;\xi(\theta)^+)+\beta(\theta;\eta(\theta)^-)\nonumber\\
+\Big(\big(\alpha(\xi(\theta)^-)-\beta(\theta;\xi(\theta)^-)\big)\nonumber\\
-\big(\alpha(\xi(\theta)^+)
-\beta(\theta;\xi(\theta)^+)\big)\Big)\xi^{\prime}(\theta).
\end{eqnarray}
Next, the starting time of a green period has the form $\xi(\theta)=kC+\theta$,
$k=1,\ldots$, and hence $\xi^{\prime}(\theta)=1$. Moreover, by Assumption 1.iv, $\alpha(\tau)$ is continuous
at $\tau=\xi(\theta)$ and hence $\alpha(\xi(\theta)^-)=\alpha(\xi(\theta)^+)$, and finally,
$\beta(\theta;\xi(\theta)^-)=0$ since $\xi(\theta)^-$ lies in a red-signal period.
All of this reduces (22) to
\begin{equation}
f^{\prime}(\theta)=\beta(\theta;\eta(\theta)^-).
\end{equation}
Consider now the case where $\xi(\theta)=u_{t}(\theta)$. Then, the lower boundary of the integral
in (16) is $u_{t}(\theta)$, and the corresponding boundary condition in (17) becomes
$-\big(\alpha(u_{t}(\theta)^+)-\beta(\theta;u_{t}(\theta)^+)\big)u_{t}^{\prime}(\theta)$ instead of
$\Big(\big(\alpha(\xi(\theta)^-)-\beta(\theta;\xi(\theta)^-)\big)
-\big(\alpha(\xi(\theta)^+)
-\beta(\theta;\xi(\theta)^+)\big)\Big)\xi^{\prime}(\theta)$.
As a result, (22) becomes
\begin{eqnarray}
f^{\prime}(\theta)=-\beta(\theta;u_{t}(\theta)^+)+\beta(\theta;\eta(\theta)^-)\nonumber\\
-\big(\alpha(u_{t}(\theta)^+)-\beta(\theta;u_{t}(\theta)^+)\big)u_{t}^{\prime}(\theta).
\end{eqnarray}
We now assert that
\begin{equation}
\big(\alpha(u_{t}(\theta)^+)-\beta(\theta;u_{t}(\theta)^+)\big)u_{t}^{\prime}(\theta)=0.
\end{equation}
Recall that $u_{t}(\theta)$ is the time a non-empty period starts at the queue.
There are three ways a non-empty period  can start: while the queue is empty, (i) $\alpha(\cdot)$ jumps up;
(ii) $\beta(\theta;\cdot)$ jump down; and (iii) $\alpha(\cdot)-\beta(\theta;\cdot)$ rises in a continuous fashion from non-positive to positive. In the first case $u_{t}(\theta)$ is a jump time of $\alpha(\cdot)$, and since the latter process is independent
 of $\theta$, $u_{t}^{\prime}(\theta)=0$. In the second case, the only way $\beta(\theta,\cdot)$
 can jump down is at the start of
 red-signal periods. In that case $u_{t}(\theta)=kC$ for some $k=1,\ldots$, and again $u_{t}^{\prime}(\theta)=0$. In the third case,
 $\alpha(u_{t}(\theta))-\beta(\theta;u_{t}(\theta))=0$. In all three cases,
 (25) is satisfied.

 Applying (25) to (24), we obtain that
 \begin{equation}
 f^{\prime}(\theta)=-\beta(\theta;u_{t}(\theta)^+)+\beta(\theta;\eta(\theta)^-).
 \end{equation}
Finally, consider the case where $t=\eta(\theta)$. Then, in Equation (21) we have that $\eta^{\prime}(\theta)=0$ as before,
and hence the derivations of Equations (23) and (26) remain unchanged.

Consider now Equation (15). For $j=1$,  if $\xi_{1}(\theta)=u_{t}(\theta)$ then Equation (26) applies
to $f_{1}^{\prime}(\theta)$. On the other hand, if $\xi_{1}(\theta)>u_{t}(\theta)$ then Equation (23) applies
to $f_{1}^{\prime}(\theta)$, but in this case (by definition of $u_{t}(\theta)$) $u_{t}(\theta)^+$ lies in a red-signal period
and hence $\beta(\theta;u_{t}(\theta)^+)=0$; implying that (26) applies as well. Thus, in any event,
$f_{1}^{\prime}(\theta)=-\beta(\theta;u_{t}(\theta)^+)+\beta(\theta;\eta_{1}(\theta)^{-})$. For every $j=2,\ldots,p$, $\xi_{j}(\theta)>u_{t}(\theta)$ and hence Equation (23) applies, namely, $f_{j}^{\prime}(\theta)=\beta(\theta;\eta_{j}(\theta)^-)$.
Therefore, by summing up all the terms in (15), Equation (12) follows.
\end{proof}
\end{proposition}

Equation (12) requires the on-line monitoring  of traffic-flow rates, and this can be done by measuring the speed of automobiles
crossing the intersection. In the special case where  the service rate alternates between 0 and $\beta_{m}$,
 $\beta(\theta;t)$ can be directly determined by the color of the traffic light.
 
 Finally, we point out that the sample function $x(\theta;t)$ evidently is  continuous in $\theta$ for every given $t$,
 and its derivative $x^{\prime}(\theta,t)$ is monotone nondecreasing in $\theta$ (see (12)). This implies that $L(\theta)$ is
 continuous as well, and the IPA derivative
 $L^{\prime}(\theta)$ is unbiased.

\section{Simulation Examples}
This section presents  simulation examples for testing the effectiveness of our regulation technique.
The traffic-light cycle is $C=1$, and the control cycle consists of 20 light cycles.
The arrival rate consists of a off/on model where, in the {\it off} stage $\alpha(t)=0$, while for each {\it on}
 stage $\alpha(t)$ is uniformly distributed in an interval
 $[(1-\zeta)\bar{\alpha},(1+\zeta)\bar{\alpha}]$; we chose its mean to be $\bar{\alpha}=4.1$, and consider
  different values of $\zeta>0$.
  The arrival rate $\alpha(t)$ varies from one
{\it on} period to the next
but retains a constant value throughout each {\it on} period. The durations of {\it off} periods and {\it on}
periods are drawn from the uniform distributions
on the intervals $[0,0.02]$ and $[0,0.063]$, respectively. The service rate $\beta(\theta;t)$ ramps up at the start of
each green-signal period at the rate of $0.62$,  until either it reaches the  saturation level of $\beta_{m}:=5.0$ or an empty period starts. In
the latter case the service rate jumps to $\beta_{m}$, and in both cases $\beta(\theta;t)$ remains at the level of $\beta_{m}$ to the end of the green-cycle period.

The set-point reference value is $r=0.3$, and the initial control variable was set to
$\theta_{1}=0.9$.
We chose $\zeta=0.3$ and $\zeta=0.1$,
respectively, and  thus,  the variance of the arrival process is  larger in the first experiment than in the second one.

Figure~3 depicts the graphs of the obtained outputs $L_{n}(\theta_{n})$ as functions of the counter $n=1,\ldots,50$,
while Figure~4 provides  the same information for $n=10,\ldots,50$ in order to highlight the effects of the variance on the asymptotic behavior of the outputs.
In both figures the blue and red graphs correspond to the respective cases  of $\zeta=0.3$ and $\zeta=0.1$.
In Figure 3 the graphs are hardly distinguishable, and both exhibit convergence to about $r=0.3$ in about 10 iterations.
In Figure 4 the differences are more evident, and the two graphs exhibit   variability about the target value of $0.3$. This
is expected in view of the variations in $L_{n}(\theta_{n})$ which are due to the random elements of the system
and the fact that it is nowhere near steady state after $20$ light cycles. However, the  respective means over the last 41 iterations, namely
the quantities $\frac{1}{41}\sum_{n=10}^{50}L_{n}(\theta_{n})$, are $0.3011$ for the case where $\zeta=0.3$, and $0.305$ for the case
where $\zeta=0.1$.

\begin{figure}
\centering
\includegraphics[width=3.6in]{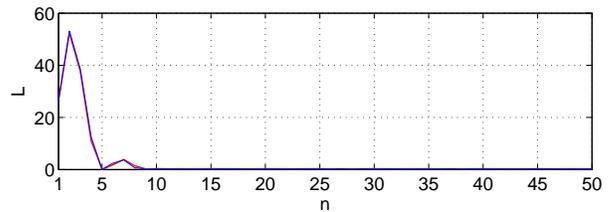}
\caption{Evolution of $L_{n}:\ n=1\;\dots,50$}
\label{Figure 3}
\end{figure}

\begin{figure}[!t]
\centering
\includegraphics[width=3.6in]{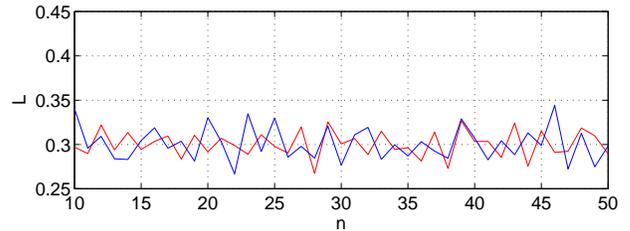}
\caption{Evolution of $L_{n}:\ n=10\;\dots,50$}
\label{Figure 1}
\end{figure}

Finally, Figure~5 shows plots of the control variable $\theta_{n}$ as functions of $n$ for the case where $\zeta=0.3$,
for two runs with the respective starting values of $\theta_{1}=0.9$ (the blue graph) and $\theta_{1}=0.1$ (the green graph). Not surprisingly, both
settle to roughly the same value ($\theta\sim 0.25$) after 10 iterations.
 We point out that the flat part of the blue curve at iterations 3-5 is due to a lower-bound guard on $\theta$
at $0.1$,
designed
to prevent extreme values which could destabilize the system.

\begin{figure}[!t]
\centering
\includegraphics[width=3.6in]{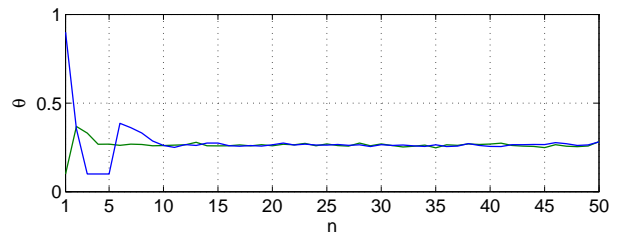}
\caption{Evolution of $\theta_{n}:\ n=1,\ldots,50$}
\label{g}
\end{figure}

\section{Conclusions and Future Work}

This paper proposes a regulation technique for congestion management in a traffic-light intersection. The technique aims at
tracking a given reference queue level at the light in the face of variable traffic patterns. It
is based on the simple idea of an integral controller with a variable gain, adjusted according to the IPA derivative of the plant function. The main theoretical result concerns a simple formula for the IPA derivative, which is computable from traffic
rates that can be measured on-line.
  Simulation results exhibit fast convergence towards the set value, and suggest the potential viability of
our approach in eventual applications.

Future research concerns  extensions of our control formulation to grids of traffic light with cross-correlated traffic patterns. On the theoretical
side, the main question is how to design  regulators for traffic-light systems with multiple controllers.
On the practical side, the main issue is how to apply the proposed technique to achieve
effective regulation under more realistic traffic conditions.

\end{document}